\documentclass[12pt, reqno]{article}
\usepackage{helvet}
\usepackage{amsfonts}
\usepackage{latexsym}
\usepackage{amsmath}
\usepackage{amsthm,amssymb,mathrsfs,amstext}
\font\msbm=msbm10

\numberwithin{equation}{section}
 
\theoremstyle{plain}
\newtheorem{Theorem}{Theorem}[section]
\newtheorem{lemma}[Theorem]{Lemma}
\newtheorem{corollary}[Theorem]{Corollary}

\newtheorem{remark}{Remark}[section]

\newtheorem{definition}{Definition}[section]

\title{Wiener Amalgam Spaces with respect to Quasi-Banach Spaces}
\author{Holger Rauhut\\
NuHAG, Faculty of Mathematics, University of Vienna\\
Nordbergstrasse 15, A-1090 Wien, Austria\\
rauhut@ma.tum.de}
\date{}

\def\mathbb#1{\hbox{\msbm{#1}}}
\newcommand{\N}{{\mathbb{N}}}
\newcommand{\R}{{\mathbb{R}}}
\newcommand{\Z}{{\mathbb{Z}}}

\newcommand{\G}{{\cal{G}}}

\newcommand{\on}{{|\!|\!|}}

\newcommand{\beq}{\begin{eqnarray}}
\newcommand{\eeq}{\end{eqnarray}}
\newcommand{\beqn}{\begin{eqnarray*}}
\newcommand{\eeqn}{\end{eqnarray*}}

\newcommand{\supp}{\operatorname{supp}}

\renewcommand{\qed}{\rule{2.5mm}{2.5mm}}
\newenvironment{Proof}{\noindent
{\bf\underline{Proof:} }}
{\hspace*{\fill}\qed\vskip1em}

\begin{document}
\maketitle
\begin{abstract}
We generalize the theory of Wiener amalgam spaces on locally compact groups 
to quasi-Banach spaces. As a main result we
provide convolution relations for such spaces. Also we weaken the
technical assumption that the global component is invariant under
right translations, which is even new for the classical Banach space case.
To illustrate our theory we discuss in detail an example on the
$ax+b$ group.
\end{abstract}

\noindent
{\bf AMS subject classification:}  46A16, 46E27, 46E30


\noindent
{\bf Key Words:}  Wiener amalgam spaces, quasi-Banach spaces, convolution
relations, doubling weights

\section{Introduction}

Wiener amalgam spaces consist of functions on a locally compact group
defined by a (quasi-)norm
that mixes, or amalgamates, a local criterion with a global criterion. 
The most general definition of Wiener amalgams so far was provided by
Feichtinger in the early 1980's in a series of papers 
\cite{Fei_Wiener,Fei_Wiener2,Fei_Homog}. We refer to \cite{Heil} 
for some historical notes and for an 
introduction for Wiener amalgams on the real line.

Wiener amalgams have proven to be a very useful tool for instance in 
time-frequency analysis \cite{Gr} (e.g. the Balian-Low theorem \cite{Heil}) and
sampling theory. Our interest in those spaces arose from coorbit space theory 
\cite{FG1,FG2,FG3,RauQB} 
which provides a group-theoretical approach to function spaces like
Besov and Triebel-Lizorkin spaces as well as modulation spaces. 

It seems that Wiener amalgams with respect to quasi-Banach spaces have not
yet been considered in full generality, except for a few results for 
Wiener amalgams on $\R^d$ in \cite{GS}. So this paper deals with
basic properties of Wiener amalgams $W(B,Y)$ 
with a quasi-Banach space $Y$ as global component
and one of the spaces $B=L^1,L^\infty$ or $M$ 
(the space of complex Radon measures) as local component. 
Moreover, we also
remove the technical assumption imposed by Feichtinger \cite{Fei_Wiener} 
that the global component $Y$ has to
be invariant under right translation. Thus,  some of our results
are even new for the classical case of Banach spaces $Y$.

One of our main achievements is a convolution relation for Wiener amalgams.
As a special case it turns out that $W(L^\infty,L^p)$ is a convolution
algebra for $0<p\leq 1$ if the underlying group is an IN group, e.g. $\R^d$.
This result is interesting since for non-discrete groups 
there are no convolution relations available for $L^p$ if $p<1$. 
The problem comes from possible $p$-integrable singularities which
are not integrable. So the integral defining the convolution $F*G$ does not
even exist for all $F \in L^p$ even if $G$ is very nice, e.g. continuous
with compact support. Of course, the local component $L^\infty$ of 
$W(L^\infty,L^p)$ prohibits such singularities. So our results indicate
that whenever treating quasi-Banach
spaces in connection with convolution then one is almost forced to use
Wiener amalgam spaces. 

To illustrate our results we also treat a class of spaces $Y$ 
on the $ax+b$ group such that
$W(L^\infty,Y)$ is right translation invariant (and thus admits convolution
relations) although $Y$ is not.


For a quasi-Banach space $(B,\|\cdot|B\|)$, we denote the quasi-norm 
of a bounded operator $T:B\to B$ by $\on T|B\on$. 
The symbol $A\asymp B$ indicates throughout
the paper that there are constants $C_1, C_2 > 0$ 
such that $C_1 A \leq B \leq C_2 A$
(independently on other expression on which $A,B$ might depend). We usually
use the symbol $C$ for a generic constant whose precise value
might be different in each occurence.


\section{Basic properties}

Let $\G$ be a locally compact group. Integration on $\G$ will always
be with respect to the left Haar measure. We denote by $L_x F(y) = F(x^{-1}y)$ and
$R_x F(y) = F(yx)$, $x,y \in \G$, the left and right translation operators.
Furthermore, let $\Delta$ be the Haar-module on $\G$.
For a Radon measure $\mu$ we denote 
$(A_x \mu)(k) \:= \mu(R_x k)$, $x \in \G$
for a continuous function $k$ with compact support. 
We may identify 
a function $F \in L^1$
with a measure $\mu_F \in M$ by $\mu_F(k) = \int F(x) k(x) dx$. 
Then it clearly holds $A_x F = \Delta(x^{-1}) R_{x^{-1}} F$.
Further, we define the 
involutions $F^\vee(x) = F(x^{-1})$, 
$F^\nabla(x) = \overline{F(x^{-1})}$, $F^*(x) = \Delta(x^{-1}) \overline{F(x^{-1})}$. 

A quasi-norm $\|\cdot\|$ 
on some linear space $Y$ is defined in the same way as a norm, with
the only difference that the triangle inequality
is replaced by $\|f+g\| \leq C(\|f\|+\|g\|)$ with some constant 
$C\geq 1$. 
It is well-known, 
see e.g. \cite[p.~20]{DL} or \cite{Kalton}, that 
there exists an equivalent quasi-norm $\|\cdot|Y\|$  
on $Y$ and an exponent $p$ with $0 < p \leq 1$ such that $\|\cdot|Y\|$
satisfies the $p$-triangle inequality, i.e., 
$\|f+g|Y\|^p \leq \|f|Y\|^p + \|g|Y\|^p$. ($C$ and $p$ are related by $C=2^{1/p}-1$.) We can choose $p=1$ 
if and only if $Y$ is a Banach space. We always
assume in the sequel that such a $p$-norm 
on $Y$ is chosen and denote it by $\|\cdot|Y\|$. If $Y$ is complete
with respect to the topology defined by the metric $d(f,g) = \|f-g|Y\|^p$
then it is called a quasi-Banach space.

Let $Y$ be a quasi-Banach space of measurable functions 
on $\G$, which contains
the characteristic function of any compact subset of $\G$. 
We assume $Y$ to be
solid, i.e., if $F \in Y$ and $G$ is measurable and satisfies $|G(x)| \leq |F(x)|$ a.e. then
also $G \in Y$ and $\|G|Y\| \leq \|F|Y\|$.

The Lebesgue spaces $L^p(\G)$, $0<p\leq \infty$ provide natural 
examples of such spaces $Y$, and the usual quasi-norm in $L^p(\G)$ is a 
$p$-norm if $0<p\leq 1$.
If $w$ is some positive measurable weight function on $\G$ then we 
further define 
$L^p_w = \{F \mbox{ measurable }, Fw \in L^p\}$ with  
$\|F|L^p_w\| := \|Fw|L^p\|$. A continuous weight $w$ is called submultiplicative
if $w(xy) \leq w(x) w(y)$ for all $x,y \in \G$. 

Now let $B$ be one of the spaces $L^\infty(\G), L^1(\G)$ or $M(\G)$, 
the space of complex Radon measures.  
Choose some relatively compact neighborhood $Q$ of $e \in \G$.
We define the control function by
\begin{equation}\label{def_control}
K(F,Q,B)(x) \,:=\, \|(L_x \chi_Q) F|B\|, \quad x\in \G,
\end{equation}
if $F$ is locally contained in $B$, in symbols $F\in B_{loc}$.
The {\bf Wiener amalgam space} $W(B,Y)$ is then defined as
\[\label{def_Wiener_space}
W(B,Y) \,:=\, W(B,Y,Q) \,:=\, \{F \in B_{loc},\, K(F,Q,B) \in Y\}
\]
with quasi-norm
\begin{equation}\label{qnormW}
\|F|W(B,Y,Q)\|\,:=\, \|K(F,Q,B)|Y\|.
\end{equation}
$B$ is called the local component and $Y$ the global component.
It follows from the solidity of $Y$ and from the quasi-norm properties 
of $\|\cdot|B\|$ and 
$\|\cdot|Y\|$ that (\ref{qnormW}) is indeed a quasi-norm. 
Since $B$ is a Banach space it is easy to see that 
also (\ref{qnormW}) is a $p$-norm 
(with $p$ being the exponent of the quasi-norm of $Y$). 
We emphasize that in general we do not require here that $Y$ is right 
translation invariant in contrast to the classical papers of 
Feichtinger \cite{Fei_Wiener,Fei_Wiener2}.

\begin{remark} The restriction of the local component $B$ to the spaces
$L^1, L^\infty$ and $M$ is done for the sake of simplicity. One can certainly
extend our considerations to more general spaces $B$, e.g. $L^p$-spaces
with $0<p\leq \infty$, compare \cite{Fei_Wiener,Heil}. 
However, convolution relations as in Section
\ref{sec_conv} will not hold any more when taking $B=L^p$ for $p<1$.
\end{remark}


Let us first make some easy observations.

\begin{lemma} We have the following continuous embeddings.
\begin{itemize}\itemsep-1pt
\item[(a)] $W(L^\infty,Y) \hookrightarrow Y$. 
\item[(b)] 
$
W(L^\infty,Y) \hookrightarrow W(L^1,Y) \hookrightarrow W(M,Y).
$
\end{itemize}
\end{lemma}
\begin{Proof} (a) Since $|F(x)| \leq \sup_{u\in U} |F(u^{-1}x)|$ for a compact neighborhood $U$ of 
$e \in \G$ the assertion follows from the solidity of $Y$.
%
%

The statement (b) follows immediately from 
$L^\infty(Q) \hookrightarrow L^1(Q) \hookrightarrow M(Q)$ 
for any compact set $Q \subset \G$. 
\end{Proof} 

Let us now investigate whether $W(B,Y,Q)$ is independent of  
$Q$ and whether it is complete. It will turn out that 
both properties are connected to the right translation 
invariance of $W(B,Y)$.
In order to clarify this we need certain 
discrete sets in $\G$ and associated sequence spaces.

\begin{definition} Let $X=(x_i)_{i\in I}$ be some 
discrete set of points in $\G$ and 
$V$ a relatively compact neighborhood of $e$ in $\G$.
\begin{itemize}\itemsep=-1pt
\item[(a)] $X$ is called $V$-dense if $\G = \bigcup_{i \in I} x_i V$.
\item[(b)] $X$ is called relatively 
separated if for all compact sets
$K \subset \G$ there exists a constant $C_K$ such that
$
\sup_{j \in I} \#\{ i\in I,\, x_iK \cap x_jK \neq \emptyset \} \leq C_K. 
$
\item[(c)] $X$ is called $V$-well-spread (or simply well-spread) if it is 
both relatively separated and $V$-dense for
some $V$.
\end{itemize}
\end{definition}
The existence of $V$-well-spread sets for arbitrarily small $V$ is proven in
\cite{Fei_Homog}. 


Given the function space $Y$, a well-spread family $X=(x_i)_{i\in I}$ and a 
relatively compact neighborhood $Q$ of $e \in \G$ 
we define the sequence space 
\begin{align}\label{def_Ydiscrete}
Y_d \,:=\, Y_d(X) \,:=\, Y_d(X,Q) \,:=&\, \{ (\lambda_i)_{i \in I}, \sum_{i\in I} 
|\lambda_i|\chi_{x_i Q} \in Y\},
\end{align}
with natural norm
$\|(\lambda_i)_{i \in I} | Y_d\| := \|\sum_{i \in I} |\lambda_i| \chi_{x_i Q}|Y\|.$
Hereby, $\chi_{x_i Q}$ denotes the characteristic function of the set $x_i Q$. 
If the quasi-norm of $Y$ is a $p$-norm,
$0 < p \leq 1$, then also $Y_d$ has a $p$-norm. 
If e.g. $Y=L^p_m$, $0 < p \leq \infty$, 
with a moderate weight $m$, then it is easily seen that 
$Y_d = \ell^p_{\tilde{m}}$ with $m(i) = \tilde{m}(x_i)$. 

Although we will not require the right translation of $Y$ in general, we state the
following easy observation in case it holds.
\begin{lemma}\label{lem_Yd_indep} 
If $Y$ is right translation invariant then the definition of \linebreak $Y_d=Y_d(X,U)$ 
does not depend on $U$.
\end{lemma}
\begin{Proof} Let $V$, $U$ be relatively compact sets with non-void interior. 
Then there exists
a finite number of points $y_j, j=1,\hdots,n$, such that $V=\cup_{j=1}^n U y_j$.
This implies 
\[
\sum_{i\in I} |\lambda_i| \chi_{x_i V} \,\leq\, \sum_{j=1}^n \sum_{i\in I} |\lambda_i|
\chi_{x_i U y_j} \,=\, \sum_{j=1}^n R_{y_j^{-1}}\left(\sum_{i\in I} |\lambda_i|\chi_{x_i U}\right).
\]
By solidity and the $p$-triangle inequality we obtain
\begin{align}
\|\sum_{i\in I} |\lambda_i|\chi_{x_iV}|Y\| \,&\leq\, 
\left( \sum_{j=1}^n \on R_{y_j^{-1}}|Y \on^p  
\|\sum_{i\in I} |\lambda_i| \chi_{x_i U}|Y\|^p\right)^{1/p}\notag\\
&=\, C \|\sum_{i \in I} |\lambda_i| \chi_{x_i U}|Y\|.\notag
\end{align}
Exchanging the roles of $V$ and $U$ shows the reverse inequality.
\end{Proof}

The following concept will also be very useful.

\begin{definition}\label{def_IBUPU} Suppose $U$ is a relatively compact
neighborhood of $e \in \G$. A collection of functions 
$\Psi = (\psi_i)_{i\in I}, \psi_i \in C_0(\G)$, is called
bounded uniform partition of unity of size $U$ 
(for short $U$-BUPU) if the 
following conditions are satisfied:
\begin{itemize}\itemsep-1pt
\item[(1)] $0 \leq \psi_i(x) \leq 1$ for all $i\in I$, $x\in \G$,
\item[(2)] $\sum_{i \in I} \psi_i(x) \equiv 1$,
\item[(3)] there 
exists a  well-spread family $(x_i)_{i\in I}$ such that 
$
\supp \psi_i \,\subset\, x_i U.
$     
\end{itemize}
\end{definition}

The construction of BUPU's with respect to arbitrary well-spread sets is standard.

We call $W(B,Y)$ right translation invariant 
if for any relatively compact neighborhood 
$Q$ of $e$ the space 
$W(B,Y,Q)$  is right 
translation invariant and the right translations 
$R_x: W(B,Y,Q) \to W(B,Y,Q)$ 
are bounded operators. (In case $B=M$ we replace $R_x$ by $A_x$ in this
definition.)

Now we are prepared to state basic properties of Wiener amalgams.

\begin{Theorem}\label{thm_basic} The following statements are equivalent:
\begin{itemize}
\item[(i)] $W(L^\infty,Y) = W(L^\infty,Y,Q)$ is independent of the choice
of the neighborhood $Q$ of $e$ (with equivalent norms for different choices). 
\item[(ii)] 
For all relatively separated sets $X$ the space 
$Y_d = Y_d(X,Q)$ is independent of the choice of the neighborhood
$Q$ of $e$ (with equivalent norms for different choices).
\item[(iii)] $W(L^\infty,Y) = W(L^\infty,Y,Q)$ is right translation invariant
(for all choices of $Q$).
\end{itemize} 
If one (and hence all) of these conditions are satisfied then also $W(B,Y)=W(B,Y,Q)$ is 
independent of the choice of $Q$. Moreover, the expression
\begin{equation}\label{WYd_norm}
\|F|W(B, Y_d)\| \,:=\, \| (\|F \psi_i|B\|)_{i\in I} |  Y_d(X)\|,
\end{equation} 
defines an equivalent quasi-norm on $W(B,Y)$, where 
$(\psi_i)_{i\in I}$ is a BUPU corresponding to the well-spread set $X$.
\end{Theorem}
\begin{Proof} We first prove that (ii) implies that (\ref{WYd_norm}) defines
an equivalent quasi-norm on $W(B,Y)$. 
Let $Q$ be a relatively compact neighborhood of $e \in \G$.
Then there exists an open set $U=U^{-1}$ with $U^2 \subset Q$.
Choose a BUPU $(\phi_i)_{i\in I}$ of size $U$.
If $x_i U \subset zQ$ then for $F \in B_{loc}$ we have
\[
\|F\phi_i|B\| \leq \|F \chi_{x_i U}|B\| \leq \|F \chi_{zQ}|B\| = K(F,Q,B)(z).
\]
This yields
\begin{equation}\label{lower_ineq}
\sum_{i\in I} \|F \phi_i|B\| \chi_{x_i U}(z) \,=\, 
\sum_{i, x_i \in zU^{-1}} \|F\phi_i|B\| \,\leq\, C K(F,Q,B)(z)
\end{equation}
since $(x_i)_{i\in I}$ is relatively separated. 
By solidity we obtain 
\[
\|(\|F\phi_i|B\|)_{i\in I}|Y_d(X,U)\| \leq C\|F|W(B,Y,Q)\|.
\]
Moreover, we have
\begin{align}
K(F,Q,B)(z) \,&=\, \|\chi_{zQ} F|B\| \,=\, \|\chi_{zQ} \sum_{i\in I} F \phi_i|B\|\notag\\
\label{upper_ineq}
&\leq\, \sum_{i, zQ \cap x_i U \neq \emptyset} \|F\phi_i|B\|
\,\leq\, \sum_{i\in I} \|F\phi_i|B\| \chi_{x_i UQ^{-1}}(z).
\end{align}
By solidity this yields
\[
\|F|W(B,Y,Q)\| \leq \|(\|F\phi_i|B\|)_{i\in I} |Y_d(X, UQ^{-1})\|.
\]
Thus, the independence of $Y_d(X,U)$ of $U$ implies that 
the norm in (\ref{WYd_norm})
is equivalent to the norm in $W(B,Y)$.  Moreover, since $Q$ was arbitrary
this shows also that $W(B,Y) = W(B,Y,Q)$ is independent of the choice of $Q$.
Specializing to $B=L^\infty$ we have thus also shown (ii) $\Longrightarrow$ (i).

As next step we prove that (iii) implies (ii). 
Let $U,V$ be relatively compact neighborhoods of $e$. Choose a neighborhood
$Q=Q^{-1}$ of $e \in G$ such that $Q^2\subset V$. 
Observe that   
\begin{align}
K(\sum_{i\in I} |\lambda_i| \chi_{x_i Q},Q)(y) \,&=\, 
\sup_{z \in yQ} \sum_{i\in I} |\lambda_i| \chi_{x_i Q}(z)
\,\leq\, \sum_{i\in I} |\lambda_i| \chi_{x_i Q^2}(y)\notag\\
&\leq\, \sum_{i\in I} |\lambda_i| \chi_{x_i V}(y).\notag
\end{align}
The right translation invariance of $W(L^\infty,Y,Q)$ together with 
Lemma \ref{lem_Yd_indep} applied to $W(L^\infty,Y)$
and the trivial inequality $|F(x)| \leq \sup_{z \in xQ}|F(z)|$ 
thus imply
\begin{align}
&\| \sum_{i\in I} |\lambda_i| \chi_{x_i U}| Y\| \,\leq\, 
\| K(\sum_{i \in I} |\lambda_i| \chi_{x_i U},Q,L^\infty)|Y\|\notag\\
&\leq\, \|K(\sum_{i \in I} | \lambda_i| \chi_{x_i Q},Q,L^\infty)|Y\|
\label{ineq_WLY}
\,\leq\, \|\sum_{i\in I} |\lambda_i|\chi_{x_i V}|Y\|.
\end{align}
Exchanging the roles of $U$ and $V$ shows the reverse inequality.


Finally, we prove (i) $\Longrightarrow$ (iii). 
Let $F \in W(L^\infty,Y)$ and $y \in \G$. We can 
find a compact neighborhood $V^{(y)}$ of $e$ 
such that $Qy \subset V^{(y)}$.
We obtain
\begin{align}
K(R_y F, Q,L^\infty)(x) \,&=\, \|(L_x \chi_Q) (R_y F)\|_\infty \,=\, \|(R_{y^{-1}} L_x \chi_Q) F\|_\infty \notag\\  
&=\, \|(L_x \chi_{Qy})F\|_\infty
\,\leq\, \|(L_x \chi_{V^{(y)}}) F\|_\infty.\notag
\end{align}
By assumption this yields together with the solidity
\begin{align}
\|R_y F|W(L^\infty,Y)\| \,&\leq\, C\|K(R_y F, Q,L^\infty)|Y\| \,\leq\, C \|K(F,V^{(y)},L^\infty)|Y \| \notag\\
&\leq\, C'(y) \|F|W(L^\infty,Y)\|.\notag
\end{align}
This concludes the proof.
\end{Proof}
 
\begin{remark}\label{rem_equivnorm} \begin{itemize}
\item[(a)]
The proof of the equivalence of the quasi-norm
in (\ref{WYd_norm}) still works (with slight changes) 
when replacing the BUPU $(\psi_i)_{i\in I}$ by the characteristic functions
$\chi_{x_i U}$. Thus, if $Y_d=Y_d(X,Q)$ is independent of the choice of $Q$
then also the expression
\[
\|(\|F\chi_{x_i Q}|B\|)_{i\in I}|Y_d\|
\]
defines an equivalent quasi-norm on $W(B,Y)$.
\item[(b)] Analyzing the proof that (ii) implies (i) one recognizes 
that it is actually
enough to require that for all neighborhoods $Q$ of $e$ there
exists some relatively separated $Q$-dense set $X$ such that 
$Y_d(X,U)$ is independent of the choice of $U$. The theorem then shows   
that $Y_d(X,U)$ is automatically 
independent of $U$ for all relatively separated sets $X$.
\end{itemize}
\end{remark} 



\begin{corollary}\label{cor_WLYd_Yd} If $W(L^\infty,Y)$ is right 
translation invariant then \linebreak
$(W(L^\infty,Y))_d=Y_d$.
\end{corollary}
\begin{Proof} This follows immediately from inequality (\ref{ineq_WLY}).
\end{Proof}

Let us now investigate the completeness of the spaces $W(B,Y)$ and $Y_d$.

\begin{lemma}\label{lem_Yd_complete} $Y_d$ is complete, and convergence in $Y_d$ implies coordinatewise
convergence.
\end{lemma}
\begin{Proof} Let $\Lambda^n = (\lambda_i^{(n)})_{i\in I}$, $n\in  \N$, be a Cauchy
sequence in $Y_d$. 
This means that the functions 
$F_n \,=\, \sum_{i\in I} \lambda_i^{(n)} \chi_{x_i U}$ form a Cauchy sequence
in $Y$. Since $Y$ is complete the limit $F = \lim_{n\in N} F_n$ exists. 
It follows from the solidity that $F$ has the form 
$F= \sum_{i\in I} \lambda_i \chi_{x_i U}$ with 
$\lambda_i = \lim_{n\to \infty} \lambda_i^{(n)}$. Clearly,
$(\lambda_i)_{i\in I} \in Y_d$ is the limit of $\Lambda^n$.
\end{Proof}

\begin{Theorem}\label{thm_complete} If $W(L^\infty,Y)$ is 
right translation invariant then $W(B,Y)$ is complete.
\end{Theorem}
\begin{Proof} Let $(\psi_i)_{i\in I}$ be some
BUPU of size $U$. By Theorem \ref{thm_basic}
$\|\cdot|W(B,Y_d)\|$ defined in (\ref{WYd_norm}) 
is an equivalent quasi-norm on $W(B,Y)$. 
Assume that $F_n$, $n\in \N$, is a Cauchy sequence of functions in $W(B,Y)$.
This implies that $(\|F_n \psi_i|B\|)_{i\in I}$ is a Cauchy sequence 
in $Y_d$ and by Lemma \ref{lem_Yd_complete} the sequence $(F_n \psi_i)_{n \in \N}$ 
is a Cauchy
sequence in $B$ for each $i \in I$. Since $B$ is complete the limit
$\lim_{n\to \infty} F_n \psi_i = F^{(i)}$ exists for each $i\in I$.
Set $F:= \sum_{i\in I} F^{(i)}$. Clearly, $\supp F^{(i)} \subset x_i U$.
Furthermore, 
\begin{align}
\|F \psi_i|B\| &\,=\, \|\sum_{j\in I} F^{(j)} \psi_i|B\| \,=\, 
\|\sum_{j: x_i U \cap x_i U} F^{(j)} \psi_i|B\|\notag\\
&\leq\, \sum_{j:x_j U \cap x_i U} \|\lim_{n\to \infty} F_n \psi_j \psi_i|B\|
\,\leq\, C \|F^{(i)}|B\|. \notag 
\end{align}
By completeness of $Y_d$, 
the sequence $(\|F^{(i)}|B\|)_{i\in I}$ is contained
in $Y_d$, and hence $F \in W(B,Y)$. 
Furthermore, we have
\[
F \,=\, \sum_{i\in I} F^{(i)} \,=\, \sum_{i\in I} \lim_{n\to \infty} F_n \psi_i
\,=\, \lim_{n\to \infty} F_n \sum_{i\in I} \psi_i = \lim_{n\to \infty} F_n.
\]
Thus, $F$ is the limit of $F_n$ in $W(B,Y)$ and hence, $W(B,Y)$ is complete.
\end{Proof}

\section{Left translation invariance}

Also the left translation invariance is an important property. 
In this section we assume that $W(L^\infty,Y)$ is right translation invariant,
so that $W(B,Y)$ is complete and independent of the choice 
of the neighborhood $Q$
according to Theorems \ref{thm_complete} and \ref{thm_basic}.

\begin{lemma}\label{lem_linfty_embed} 
If $W(L^\infty,Y)$ is left translation invariant then
$Y_d$ is continuously embedded into $\ell^\infty_{1/r}$ with 
$r(i):= \on L_{x_i^{-1}}|W(L^\infty,Y)\on$.
\end{lemma} 
\begin{Proof} Let $U$ be some compact neighborhood of $e$ and  
$(\lambda_i)_{i\in I} \in Y_d$.
With $C:= \|\chi_U|W(L^\infty,Y)\|$ we obtain
by Corollary \ref{cor_WLYd_Yd} and solidity
\begin{align}
C |\lambda_i| &=\, |\lambda_i| \|\chi_U|W(L^\infty,Y)\| \,=\, |\lambda_i|
\|L_{x_i^{-1}} \chi_{x_i U}|W(L^\infty,Y)\|  \notag\\
&\leq \,\on L_{x_i^{-1}}|W(L^\infty,Y)\on
\||\lambda_i| \chi_{x_i U}|W(L^\infty,Y)\|\notag\\ 
&\leq\, r(i) \|\sum_{j\in I} |\lambda_j| \chi_{x_j U}|W(L^\infty,Y)\|
\,\leq\, r(i) \|(\lambda_i)_{i\in I}|Y_d\|.\notag
\end{align}
This completes the proof.
\end{Proof}


\begin{lemma}\label{lem_Linfty_embed}
If $W(L^\infty,Y)$ is left translation invariant then
$W(L^\infty,Y)$ is continuously embedded into $L^\infty_{1/r}$, where 
$r(x):= \on L_{x^{-1}}|W(L^\infty,Y)\on$. 
\end{lemma} 
\begin{Proof}
By Theorem \ref{thm_basic} 
$Y_d=Y_d(X,Q)$ 
is independent of the choice of $Q$ and 
the quasi-norm $\|\cdot|W(L^\infty,Y_d)\|$ defined in (\ref{WYd_norm})
is equivalent to the quasi-norm of $W(L^\infty,Y)$. Since $Y_d$ is continuously
embedded into $\ell^\infty_{1/r}$ by Lemma 
\ref{lem_linfty_embed} and $(L^\infty_{1/r})_d = \ell^\infty_{1/r}$ we obtain
\begin{align}
C_1\|F|W(L^\infty,L^\infty_{1/r})\| \,&\leq\, 
\|F|W(L^\infty,\ell^\infty_r)\| \,\leq\, \|F|W(L^\infty,Y_d)\|\notag\\
&\leq\, C_2 \|F|W(L^\infty,Y)\|
\end{align}
for all $F \in W(L^\infty,Y)$. Further, it is easy to see that 
$W(L^\infty,L^\infty_{1/r}) = L^\infty_{1/r}$.
\end{Proof}

In some cases one has translation invariant spaces $Y$. Then
we have the following estimates of the norm of the left translation operators
in $W(L^\infty,Y)$.

\begin{lemma}\label{lem_left_trans} If $Y$ is left translation invariant then $W(B,Y)$
is left translation invariant and 
$
\on L_y |W(B,Y)\on \leq \on L_y | Y \on.
$
\end{lemma}
\begin{Proof} We have
\begin{align}
K(L_y F,Q,B)(x) \,&=\, \|(L_x \chi_Q) (L_y F)|B\|
\,=\, \|(L_{y^{-1}x} \chi_Q) F|B\|\notag\\
&=\, (L_y K(F,Q,B))(x).\notag
\end{align}
This yields
\[
\|L_y F|W(B,Y)\| = \|L_y K(F,Q,B)|Y\| \leq \on L_y|Y\on \|F|W(B,Y)\|,
\]
and the proof is completed.
\end{Proof}

\section{Conditions ensuring translation invariance}

Given a concrete space $Y$, according to the previous results, 
there is the need to check whether $W(L^\infty,Y)$ 
is right translation invariant. Moreover, we will see later that also
the right translation invariance of $W(M,Y)$ is important in order to have convolution
relations.


\begin{lemma}\label{lem_WMY} If $W(L^\infty,Y)$ is right translation 
invariant then also $W(M,Y)$ is right translation invariant.
\end{lemma}
\begin{Proof}
Let $\mu \in W(M,Y)$, $y \in \G$ and $Q$ be a compact neighborhood of
$e$. Then there exist a finite number of points $y_k, k=1,\hdots,n$, such that 
such that $Qy^{-1} \subset \bigcup_{k=1}^n y_k Q$. We obtain for the control function
\begin{align}
K(A_y \mu,Q,M)(x) \,&=\, \|(L_x \chi_Q) A_y \mu|M\|
\,=\, |\mu|(R_y L_x \chi_Q) \,=\, |\mu|(L_x \chi_{Qy^{-1}})\notag\\
 &\leq\,
\sum_{k=1}^n |\mu|(L_x \chi_{y_k Q}) \,=\, \sum_{k=1}^n R_{y_k} K(\mu,Q,M)(x).\notag
\end{align}
By solidity, the $p$-triangle inequality and independence of $W(M,Y,Q)$ of
the choice of $Q$ we get
\begin{align}
&\|A_y \mu|W(M,Y)\|^p \leq\, \| \sum_{k=1}^n R_{y_k} K(\mu,Q,M)|Y\|^p\notag\\
&\leq\, \sum_{k=1}^n \|R_{y_k} K(\mu,Q,M)|W(L^\infty,Y)\|^p\notag\\
&\leq\, \sum_{k=1}^n \on R_{y_k}|W(L^\infty,Y)\on^p \|K(\mu,Q,M)|W(L^\infty,Y)\|^p\notag\\
&\leq\, \sum_{k=1}^n \on R_{y_k}|W(L^\infty,Y)\on^p \|K(\mu,Q^2,M)|Y\|^p\notag\\
&\leq\, C \sum_{k=1}^n \on R_{y_k}|W(L^\infty,Y)\on^p \|\mu|W(M,Y)\|^p.\notag
\end{align}
This concludes the proof.
\end{Proof}


Let us give another criterion for the right translation invariance of $W(B,Y)$.

\begin{corollary} If $Y$ is right translation invariant then also 
$W(B,Y) = W(B,Y,Q)$ is right 
translation invariant and independent of $Q$.
\end{corollary}
\begin{Proof} By Lemma (\ref{lem_Yd_indep}) $Y_d=Y_d(X,U)$ is 
independent of $U$. Thus,
Theorem \ref{thm_basic} implies that $W(B,Y)=W(B,Y,Q)$ is independent of $Q$ and
$W(L^\infty,Y)$ is right translation invariant. Lemma \ref{lem_WMY} implies that
also $W(M,Y)$ is right translation invariant. Clearly, $W(L^1,Y)$ is a subspace
of $W(M,Y)$ that is right translation invariant if $W(M,Y)$ is right translation invariant.
Thus, we proved the assertion for all admissible choices $B=L^\infty,L^1,M$.
\end{Proof}

Recall that $\G$ is called an IN group if there exists a compact neighborhood of $e$
such that $xQ = Qx$ for all $x \in \G$. 

\begin{lemma}\label{lem_IN_trans} 
Let $\G$ be an IN group and assume $Y$ to be right translation
invariant. Then it holds $\on R_y|W(L^\infty,Y)\on \leq \on R_y|Y\on$ and
$\on A_y|W(M,Y)\on \leq \on R_y|Y\on$.
\end{lemma}
\begin{Proof} Choose $Q$ to be a compact invariant neighborhood of $e$, i.e.,
$yQ=Qy$ for all $y \in \G$. This yields
\begin{align}
K(R_y F,Q,L^\infty)(x) \,&=\, \|(L_x \chi_Q)R_y F\|_\infty
\,=\, \|(L_x \chi_{Qy}) F\|_\infty \,=\, \|(L_{x} \chi_{yQ}) F\|_\infty\notag\\
 &=\, \|(L_{xy} Q)F\|_\infty \,=\, K(F,Q,L^\infty)(xy).\notag
\end{align}
and thus,
\[
\|R_yF|W(L^\infty,Y)\| \,=\, \|R_y K(F,Q,L^\infty)|Y\| \,\leq\, \on R_y |Y\on \|F|W(L^\infty,Y)\|.
\]
The proof for $B=M$ 
is similar. 
%
\end{Proof}

We remark that $Y$ does not necessarily need to be translation invariant in 
order $W(L^\infty,Y)$ to be
translation invariant, see Section \ref{sec_ex}.
The following criterions allow to check left or right translation invariance 
of $W(L^\infty,Y)$
without using translation invariance of $Y$.

\begin{lemma}\label{lem_left_invariant} Let $U$ be some compact neighborhood of $e \in \G$.
Let $X=(x_i)_{i\in I}$ be some well-spread set in $\G$. Denote by $x^{-1}X$, $x \in \G$, the  
well-spread set $(x^{-1} x_i)_{i\in I}$. If there is a function $k(x)$ such that
\[
\|(\lambda_i)_{i\in I}|Y_d(x^{-1}X,U)\| \,\leq\, k(x) \|(\lambda_i)_{i\in I}|Y_d(X,U)\| 
\]
for all $(\lambda_i)_{i\in I} \in Y_d(X)$ then $W(B,Y)$ is left-translation invariant with
$$
\on L_{x} |W(B,Y)\on \, \leq\, C k(x). 
$$
\end{lemma}
\begin{Proof} Let $(\psi)_{i\in I}$ be some BUPU corresponding to $X$.
Since (\ref{WYd_norm}) defines an equivalent norm
on $W(B,Y)$ we obtain 
\begin{align}
\|L_x F|W(B,Y)\| \,&\leq\, C\| (\|(L_x F) \psi_i|B\|_{i\in I}|Y_d(X,U)\|\notag\\
&\leq\, C\| (\|F (L_{x^{-1}}\psi_i)|B\|)_{i\in I}|Y_d(X,U)\|\notag\\
\,&\leq\, C k(x) \,\| (\|F (L_{x^{-1}}\psi_i)\|)_{i\in I}|Y_d(x^{-1}X,U)\|.\notag
\end{align}
The system $(L_{x^{-1}} \psi_i)_{i\in I}$ is a BUPU corresponding to the well-spread set
$x^{-1} X$. 
Thus, using once more the equivalence of the norm (\ref{WYd_norm}) 
with the norm in $W(B,Y)$ we obtain
$\|L_x F|W(B,Y)\| \,\leq\, C'k(x) \|F|W(B,Y)\|.$
\end{Proof}
\begin{remark}\label{rem_trans_norm} 
If $Y_d(X,U)$ is independent of the choice of the neighborhood $U$ 
then we already know from Theorem \ref{thm_basic} that $W(L^\infty,Y)$ is right 
translation invariant.
If $h(x)$ is a function such that
\[
\|(\lambda_i)_{i\in I}|Y_d(X,Ux)\| \,\leq\, h(x) \|(\lambda_i)_{i\in I}|Y_d(X,U)\|
\]
for all $(\lambda_i)_{i\in I} \in Y_d(X)$ then a similar argument as in the previous proof shows that
 $$\on R_x|W(L^\infty,Y)\on \, \leq\, C h(x).$$
\end{remark}

\section{Convolution relations}
\label{sec_conv}


Let us now prove the main results of this article concerning
convolution relations of Wiener amalgams with quasi-Banach spaces
as global component (compare \cite{FG1,FG2} 
for the classical case of Banach spaces).

\begin{Theorem}\label{thm_conv} Let $0 < p \leq 1$ be 
such that the quasi-norm of 
$Y$ satisfies the $p$-triangle inequality and 
assume that $W(L^\infty,Y)$ is right translation 
invariant.
\begin{itemize}
\item[(a)] Set $w(x):= \on A_x|W(M,Y) \on$. Then we have 
\[
W(M,Y) * W(L^\infty,L^p_w) \hookrightarrow W(L^\infty,Y)
\]
with corresponding estimate for the quasi-norms.
\item[(b)] Set $v(x):= \Delta(x^{-1}) \on R_{x^{-1}} | W(L^\infty,Y)\on$. Then we have
\[
W(L^\infty,Y) * W(L^\infty,L^p_v) \hookrightarrow W(L^\infty,Y)
\]
with corresponding estimate for the quasi-norms. 
\end{itemize} 
\end{Theorem}
\begin{Proof} (a) It follows from Theorem \ref{thm_basic} 
that any $G \in W(L^\infty,L^p_w)$
has a decomposition $G = \sum_{i \in I} L_{x_i} G_i$ with $G_i \in L^\infty$,
$\supp G_i \subset Q= Q^{-1}$ for some compact $Q$ and $\sum_{i\in I} \|G_i\|^p_\infty w(x_i)^p 
\leq C \|G|W(L^\infty,L^p_w)\|^p < \infty$.

For $\mu \in W(M,Y)$ we estimate the control function of $\mu * (L_{x_i} G_i)$ by
\begin{align}
&K(\mu * (L_{x_i} G_i),Q,L^\infty)(x) 
\,=\, \sup_{z \in xQ} |\mu * (L_{x_i} G_i)(z)|\notag\\
&=\, \sup_{z \in xQ} |\int (L_y L_{x_i} G_i)(z) d\mu(y)|
\,\leq\, \|G_i\|_\infty \sup_{q \in Q} \int L_{y x_i} \chi_Q(xq)d|\mu|(y)\notag\\
&\leq\, \|G_i\|_\infty \int \chi_{Q^2}((yx_i)^{-1}x) d|\mu|(y)
\,=\, \|G_i\|_\infty \int \chi_{Q^2} (x^{-1} yx_i)d|\mu|(y)\notag\\
&= \|G_i\|_\infty \int R_{x_i} L_x \chi_{Q^2}(y) d|\mu|(y)
\,=\, \|G_i\|_\infty \|(L_x \chi_{Q^2}) (A_{x_i} \mu)|M\| \notag\\
&=\, \|G_i\|_\infty K(A_{x_i} \mu, Q^2, M)(x).\notag
\end{align}
Thus, we have
\begin{align}
\|\mu * L_{x_i} G_i|W(L^\infty,Y)\|
\,&\leq\, \|G_i\|_\infty \|K(A_{x_i} \mu, Q^2, M)|Y\|\notag\\
&\leq\, C\|G_i\|_\infty \|A_{x_i} \mu |W(M,Y)\|.\notag
\end{align}
Pasting the pieces together yields
\begin{align}
&\|\mu * G|W(L^\infty,Y)\|^p \,=\, \|\sum_{i\in I} \mu * L_{x_i} G_i|W(L^\infty,Y)\|^p\notag\\
&\leq\, \sum_{i\in I} \|\mu * L_{x_i} G_i|W(L^\infty,Y)\|^p
\label{Ax}
\,\leq\, C \sum_{i\in I} \|G_i\|^p_\infty \|A_{x_i} \mu|W(M,Y)\|^p\\
&\leq\,  C \sum_{i\in I} \|G_i\|^p_\infty 
\on A_{x_i} |W(M,Y)\on^p \|\mu|W(M,Y)\|^p\notag\\
&\leq\, C \|\mu|W(M,Y)\|^p \|G|W(L^\infty,L^p_w)\|^p.\notag
\end{align}
(b) Since $W(L^\infty,Y) \subset W(M,Y)$ all the computations done in (a) are still valid. We only
have to replace $\|A_{x_i} \mu|W(M,Y)\|$ by
$\|A_{x_i} \mu|W(L^\infty,Y)\| = \Delta(x_i^{-1}) \|R_{x^{-1}} \mu|W(L^\infty,Y)\|$ in (\ref{Ax})
to deduce (b).
\end{Proof}


\begin{Theorem}\label{thm_convYvee} Assume $Y$ is such that $W(L^\infty,Y)$ is left and right 
translation invariant. Set
$v(x):= \on L_{x^{-1}}|W(L^\infty,Y)\on$. Then
\[
W(L^\infty,L^p_v) * W(L^\infty,Y^\vee)^\vee \hookrightarrow W(L^\infty,Y).
\] 
\end{Theorem}
\begin{Proof} 
%
Let $F \in W(L^\infty, L^p_v)$ and $G \in W(L^\infty,Y)$. Similarly as in the
proof of Theorem \ref{thm_conv} we may write $F = \sum_{i\in I} L_{x_i} F_i$
with $\supp F_i \subset Q=Q^{-1}$ (compact) and 
$\sum_{i\in I} \|F_i\|^p_\infty v(x_i)^p\leq C\|F|W(L^\infty,L^p_v)\|$.
We obtain
\begin{align}
&K(F_i*G, Q, L^\infty)(x) \,=\, \sup_{z \in xQ} |F_i*G(z)|
\,\leq\, \sup_{z \in x Q} |\int_{x_i Q} F_i(y) L_{y} G(z) dy|\notag\\
&\leq\, \|F_i\|_\infty \sup_{q \in Q} \int \chi_Q(y) |(R_q G)(y^{-1}x)|dy
\,\leq\, C\|F_i\|_\infty \int \chi_{Q^2} (y) |G^\vee(x^{-1} y)| dy\notag\\
&\leq\, C\|F_i\|_\infty \int L_{x^{-1}}\chi_{Q^2} (y) |G^\vee(y)| dy
\,\leq\, C'\|F_i\|_\infty K(G^\vee, Q^2,L^\infty)(x^{-1})\notag.
\end{align}
This yields
\begin{align}
\|F_i * G|W(L^\infty,Y)\| 
\,&\leq\, C \|F_i\|_\infty \|K(G^\vee,Q^2,L^\infty)^\vee|Y\|\notag\\
&\leq\, C\|F_i\|_\infty \|G|W(L^\infty,Y^\vee)^\vee\|.\notag
\end{align}
Pasting the pieces together we get
\begin{align}
\|F*G|W(L^\infty,Y)\|^p \,&=\, \|\sum_{i\in I} (L_{x_i} F_i) * G|W(L^\infty,Y)\|^p\notag\\
&\leq\, \sum_{i\in I} \|L_{x_i} (F_i * G)|W(L^\infty,Y)\|^p\notag\\
&\leq\, C \sum_{i\in I} \on L_{x_i}|W(L^\infty,Y)\on^p \|
\,\|F_i\|_\infty^p \|G|W(L^\infty,Y^\vee)^\vee\|^p\notag\\
&\leq\, C' \|F|W(L^\infty,L^p_v)\|^p\,\|G|W(L^\infty,Y^\vee)^\vee\|^p.\notag
\end{align}
This concludes the proof.
\end{Proof}

From the previous theorem we see that the involution ${ }^\vee$ 
has some relevance.
In the case of IN groups we have the following result.

\begin{lemma}\label{lem_IN_inv}
If $\G$ be an IN group 
then $W(L^\infty,Y^\vee)^\vee = W(L^\infty,Y)$ with equivalent norms.
\end{lemma}
\begin{Proof} Let $Q$ be an invariant compact neighborhood of $e$. Then
also $Q^{-1}$ is invariant. For the control function we obtain
\begin{align}
&K(F^\vee,Q,L^\infty)(x) \,=\, \|(L_x \chi_Q) F^\vee\|_\infty
\,=\, \|(L_x \chi_Q)^\vee F\|_\infty
\,=\, \|(R_x \chi_{Q^{-1}}) F\|_\infty\notag\\
&=\, \|\chi_{Q^{-1}x^{-1}} F\|_\infty
\,=\, \|\chi_{x^{-1} Q^{-1}} F \|_\infty
\,=\, K(F,Q^{-1},L^\infty)(x^{-1}).\notag
\end{align}
This shows the claim.
\end{Proof}

Theorem \ref{thm_convYvee} implies a convolution relation for Wiener
amalgam spaces with respect to weighted $L^p$-spaces.

\begin{corollary}\label{cor_conv_Lp} 
Let $w$ be a submultiplicative weight and $0<p\leq 1$. Then it holds
\[
W(L^\infty,L^p_w) * W(L^\infty,L^p_{w^*})^\vee \hookrightarrow W(L^\infty,L^p_w).
\] 
In particular, if $\G$ is an IN-group then 
$W(L^\infty,L^p_w) * W(L^\infty,L^p_w) \hookrightarrow W(L^\infty,L^p_w)$ with corresponding
quasi-norm estimate.
\end{corollary}
\begin{Proof} The first assertion is a direct consequence of Theorem \ref{thm_convYvee} and the second
assertion follows then from Lemma \ref{lem_IN_inv}.
\end{Proof}
In particular, if $\G$ is an IN group then 
$W(L^\infty,L^p_w)$, $0\leq p \leq 1$, is a quasi-Banach algebra
under convolution. Since commutative groups are clearly IN groups this 
result applies in particular to Wiener amalgams on $\G=\R^d$.
Moreover, if $\G$ is discrete then we recover
the well-known relation \linebreak
$\ell^p_w(\G) * \ell^p_w(\G) \hookrightarrow \ell^p_w(\G)$, $0<p\leq 1$.

\section{An example on the $ax+b$ group}
\label{sec_ex}

In this section we provide an example of a non-translation invariant
space $Y$ such that $W(L^\infty,Y)$ is right translation invariant.
We consider the $n$-dimensional $ax+b$ group
$\G = \R^n \rtimes \R_+^*$ 
where $\R_+^*$ denotes the multiplicative group of positive real numbers.
The group law in $\G$ reads
$(x,a)\cdot(y,b) = (x+ay,ab)$.
The $ax+b$ group has left Haar-measure
\beqn
\int_\G f(x) dx \,=\, \int_{\R^n} \int_0^\infty f(x,a) 
\frac{da}{a^{n+1}} dx
\eeqn
and modular function $\Delta(x,a,A) = a^{-n}$. The $ax+b$ group plays
an important role in wavelet analysis and the theory of Besov
and Triebel-Lizorkin spaces.

Let $0 < p,q \leq \infty$. With some positive measurable weight 
function $v$ on $\G$ 
we define the mixed norm space $L^{p,q}(v)$
on $\G$ as the collection of measurable functions whose quasi-norm
\[
\|F|L^{p,q}(v)\| \,:=\, \left(\int_{0}^\infty \left( \int_{\R^n} |F(x,a)|^p 
v(x,a) dx\right)^{q/p} \frac{da}{a^{n+1}}\right)^{1/q} 
\]
is finite (with obvious modification in the cases $p= \infty$ or $q=\infty$).
This quasi-norm is actually an $r$-norm where $r:= \min\{1,p,q\}$.
If $v\equiv 1$ we write $L^{p,q}$. If $p=q$ then clearly,
$L^{p,p} = L^p(\G)$. It is easy to see by an integral 
transformation that $L^{p,q}$ is invariant under left and 
right translations. We remark that for reasons to become clear later 
$v$ is treated as a measure here,
so if $v$ does not vanish on a set of positive measure
then $L^{\infty,\infty}(v) = L^{\infty}(\G)$.

With a similar argument as in \cite[Proposition 2.4]{Heil},
see also \cite{Fei_Gewicht}, 
one shows (using the right translation 
invariance of the unweighted $L^{p,q}$ space) 
that $L^{p,q}(v)$, $0<p,q<\infty$, 
is right translation invariant if and only if
\begin{equation}\label{cond_right_invariant}
v((x,a)\cdot(y,b)) \,\leq\, v(x,a) w(y,b)
\end{equation}   
for some submultiplicative function $w$ (possibly depending on $p,q$). 
Now assume that $v(x,a)$ is a function
of $x$ only. Then condition (\ref{cond_right_invariant}) means
that the quotient
\begin{equation}\label{quot}
\frac{v((x,a)(y,b))}{v(x,a)} \,=\, \frac{v(x+ay)}{v(x)} 
\end{equation}
is bounded by a submultiplicative function $w$ of $y$ only. However,
since the right hand side depends also on $a \in (0,\infty)$ this can be
satisfied only in special cases (e.g. if $v$ is bounded from above and below). 
In particular, the typical choice
$v_s(x,a) = v_s(x) = (1+|x|)^s$, $s\in \R$, does not satisfy 
(\ref{cond_right_invariant}) for any submultiplicative weight $w$ on $\G$ 
if $s\neq 0$ (although it is even submultiplicative as function on $\R^n$ if 
$s \geq 0$.) In particular, $L^{p,q}(v)$ is not right 
translation invariant for
many non-trivial choices of $v$.

In the following we introduce a class of weight functions $v$
for which \linebreak $W(L^\infty,L^{p,q}(v))$ is right translation invariant.
This class, however, contains weights $v$ that do not satisfy 
\ref{cond_right_invariant}. i.e., $L^{p,q}(v)$ is not right translation
invariant, in general. 

Let $B(x,r)$ denote the ball in $\R^n$ of radius $r$ centered at $x \in \R^n$.
A positive measurable weight function $v$ on $\R^n$ is said to satisfy the 
doubling condition if there exists a constant $C$ such that
\begin{equation}\label{cond_double}
\int_{B(x,2r)} v(y) dy \,\leq\, C \int_{B(x,r)} v(y) dy
\end{equation}
for all $x\in \R^n$ and $r\in (0,\infty)$. This condition is equivalent to the 
existence of constants $c,\alpha$ such that 
\begin{equation}\label{equiv_double}
\int_{B(x,tr)} v(y) dy \,\leq\, c t^\alpha \int_{B(x,r)} v(y) dy
\quad\mbox{ for all }
x \in \R^n, r \in (0,\infty), t\geq 1.
\end{equation}
For instance weights in the 
Muckenhoupt classes $A_p$, $p>1$, satisfy the doubling 
condition \cite{Muck1}. 
A typical example of a weight in $A_\infty = \cup_{p>1} A_p$ 
is $v^{(s)}(x)=|x|^s$, $s > -1$. So doubling weights may have zeros or poles.
A further example of a doubling weight is $v_s(x) = (1+|x|)^s$, $s \in \R$.
We remark that this weight function is not contained in $A_\infty$ if 
$s \leq -1$. For another construction of
a doubling weight which is not contained in $A_\infty$ we refer 
to \cite{Muck1}.

We extend a doubling weight $v$ on $\R^n$ to $\G = \R^n \rtimes \R_+^*$ by
setting $v(x,t) = v(x)$ for $(x,t) \in \G$. Let $L^{p,q}(v)$ be the associated
mixed norm space as defined above. We will use 
Theorem \ref{thm_basic} to prove that $W(L^\infty,L^{p,q}(v))$ is
right translation invariant. In particular, let us study the associated 
sequence space $(L^{p,q}(v))_d$.

\begin{lemma} Let $0< p<\infty$, $0 < q \leq \infty$ and $v$ be a weight 
function on $\R^n$ 
Let 
$X=(x_{k,j},a_j)_{(k,j) \in I:=\Z^n \times \Z}$ be some
well-spread set in $\G=\R^n \rtimes \R_+^*$. 
If $v$ satisfies the doubling condition (\ref{cond_double}) then 
$(L^{p,q}(v))_d = (L^{p,q}(v))_d(X,U)$ is
independent of the choice of the neighborhood $U$ of $e$ in $\G$, 
and an equivalent
norm on $(L^{p,q}(v))_d(X)$ is given by
\[
\|(\lambda_i)_{i\in I}|\ell^{p,q}(\tilde{v})\| \,=\, 
 \left( \sum_{j \in \Z} \left(\sum_{k \in \Z^n} |\lambda_{k,j}|^p 
\tilde{v}_{k,j}\right)^{q/p} a_j^{-n}\right)^{1/q}
\]
where $\tilde{v}_{k,j} = \int_{B(x_{k,j},a_j)} v(y) dy$ (with the usual
modification for $q=\infty$).

Moreover, $W(L^\infty,L^{p,q}(v))$ is right translation invariant if and only
if $v$ satisfies the doubling condition.
\end{lemma}
\begin{Proof} It satisfies to show the assertion for neighborhoods
of the form $U(r,\beta) = B(0,r) \times (\beta^{-1},\beta) \subset \G$ with 
$r \in (0,\infty)$ and $\beta \in (1,\infty)$ since for an arbitrary compact 
neighborhood $U$ of $e = (0,1) \in \G$ we can 
find $r_1,r_2,\beta_1,\beta_2$ such
that $U(r_1,\beta_1) \subset U \subset U(r_2,\beta_2)$.
Observe that 
\[
(x,a) U(r,\beta) \,=\, B(x,a r) \times (a\beta^{-1},a\beta).
\]
Using the relative separatedness of $X$ we obtain for $0<q < \infty$
\begin{align}
&\|(\lambda_i)_{i\in I}|(L^{p,q}(v))_d(X,U(r,\beta))\|\notag\\
&=\, \left(\int_0^\infty \left( \int_{\R^n} 
\sum_{j\in \Z} \sum_{k\in \Z^n} |\lambda_{k,j}|^p \chi_{B(x_{k,j},a_j r)}(y)\chi_{(a_j\beta^{-1},a_j\beta)}(a) v(y) dy\right)^{q/p} \frac{da}{a^{n+1}}\right)^{1/q}\notag\\
&\asymp\,\left(\sum_{j\in \Z} \left(\sum_{k\in \Z^n} |\lambda_{k,j}|^p \int_{B(x_{k,j},a_jr)} v(y) dy\right)^{q/p} \int_{a_j\beta^{-1}}^{a_j\beta} \frac{da}{a^{n+1}}\right)^{1/q}\notag\\
&\asymp\, \left( \sum_{j \in \Z} \left(\sum_{k \in \Z^n} |\lambda_{k,j}|^p 
\int_{B(x_{k,j},a_j r)} v(y) dy\right)^{q/p} a_j^{-n}\right)^{1/q}.\notag  
\end{align}
The computation for $q=\infty$ is similar.
Thus, $(L^{p,q}(v))_d(X,U(r,\beta))$ is independent of $r$ and $\beta$ if and 
only if for all $r,s \in (0,\infty)$ there
exist constants $C_1(r,s), C_2(r,s)>0$ such that
\begin{equation}\label{ineq1}
C_1(r,s) \int_{B(x_{k,j},a_j r)} v(y) dy
\,\leq\, \int_{B(x_{k,j},a_j s)} v(y) dy
\,\leq\, C_2(r,s)\int_{B(x_{k,j},a_j r)} v(y) dy
\end{equation}
for all $(k,j) \in \Z^n \times \Z$. Let us assume without loss of generality
that $r\leq s$. Then the first inequality is clear. Moreover, by the doubling 
condition, resp. its equivalent form (\ref{equiv_double}) we have
\[
\int_{B(x_{k,j},a_j s)} v(y) dy
\,\leq\, c (s/r)^\alpha \int_{B(x_{k,j},a_j r)} v(y) dy.
\] 
So (\ref{ineq1}) is satisfied with $C_1(r,s) = 1$ and 
$C_2(r,s) = c (s/r)^\alpha$.

Since we may choose relatively separated sets of the form $(x_{j,k},a_j)$ 
of arbitrarily small
density -- e.g.~$(ab^{-j}k,b^{-j})_{k\in \Z^n,j\in \Z}$ 
with small $a>0, b>1$ --
$W(L^\infty,L^{p,q}(v))$ is right translation invariant by 
Theorem \ref{thm_basic} and Remark \ref{rem_equivnorm}(b) 
if $v$ is doubling. Conversely, if 
$W(L^\infty,L^{p,q}(v))$ is right translation invariant then (\ref{ineq1})
must hold for any choice of the relatively separated set $X = (x_{j,k},a_j)$ 
by Theorem \ref{thm_basic}. In particular, choosing $s=2,r=1$ in (\ref{ineq1})
we obtain
\[
\int_{B(x,2a)} v(y) dy \,\leq\, C_2 \int_{B(x,a)} v(y) dy
\]
for all $x\in \R^n$, $a \in (0,\infty)$, which clearly is the 
doubling condition.
\end{Proof}

Since $L^{\infty,q}(v) = L^{\infty,q}$ the analogue of the Theorem for 
$p=\infty$ is trivial. It seems that in general $W(L^\infty,L^{p,q}(v))$
is not left invariant.

In order to state the convolution relation in Theorem \ref{thm_conv} for 
our case we
estimate the norm of the right translation operators on 
$W(L^\infty,L^{p,q}(v))$ using Remark \ref{rem_trans_norm}.
Let $U=U(r,\beta)$, $r>0,\beta>1$, be a neighborhood of $e = (0,1)$ as in the previous
proof. For $(x,a),(y,b) \in \G$ we obtain
\begin{align}
&(x,a)\cdot U(r,\beta)\cdot (y,b) \,=\, \left(B(x,ar)\times a(\beta^{-1},\beta)\right)\cdot(y,b)
\notag\\ 
&=\,\{(z+sy,sb):\, z \in B(x,ar), s \in ab(\beta^{-1},\beta)\}\notag\\
&\subset\, \bigcup_{s\in a(\beta^{-1},\beta)} B(x+sy,ar) \times ab(\beta^{-1},\beta)
\,\subset\, B(x,a(\beta|y|+r)) \times ab(\beta^{-1},\beta).\notag
\end{align}
Let $X=(x_{k,j},a_j)$ be a relatively separated set in $\G$.
Proceeding as in the previous proof we deduce
\begin{align}
&\|(\lambda_i)_{i\in I}|(L^{p,q}(v))_d(X,U(r,\beta)\cdot(y,b))\|\notag\\
&\leq\, C \left(\sum_{j\in \Z} \left(\sum_{k\in \Z^n} |\lambda_{k,j}|^p 
\int_{B(x_{k,j},a_jr (\frac{\beta}{r}|y| + 1))} v(y) dy\right)^{q/p} 
\int_{a_jb\beta^{-1}}^{a_jb\beta} \frac{da}{a^{n+1}}\right)^{1/q}\notag\\
&\leq\, C  \left(\sum_{j\in \Z} \left(\sum_{k\in \Z^n} |\lambda_{k,j}|^p 
(\frac{\beta}{r}|y|+1)^\alpha \int_{B(x_{k,j},a_j r)} v(y) dy\right)^{q/p}
b^{-n} a_j^{-n}\right)^{1/q}\notag\\
&\leq\, C (1+|y|)^{\alpha/p} b^{-n/q} 
\|(\lambda_i)_{i\in I}|(L^{p,q}(v))_d(X,U(r,\beta))\|,\notag
\end{align}
where $\alpha$ is the exponent from (\ref{equiv_double}). By Remark 
\ref{rem_trans_norm} we conclude that
\[
\on R_{(y,b)}|W(L^\infty,L^{p,q}(v))\on \,\leq\, C (1+|y|)^{\alpha/p}b^{-n/q},
\]
and since $(y,b)^{-1} = (-b^{-1}y,b^{-1})$ we have
\[
\Delta((y,b)^{-1}) \on R_{(y,b)^{-1}}|W(L^\infty,L^{p,q}(v))\on
\,\leq\, C b^{n(1+1/q)} (1+b^{-1}|y|)^{\alpha/p}.
\] 
Set $w(y,b):= b^{n(1+1/q)}(1+b^{-1}|y|)^{\alpha/p}$ and $r:=\min\{1,p,q\}$. 
Then Theorem \ref{thm_conv} tells us that
\[
W(L^\infty,L^{p,q}(v)) * W(L^\infty,L^r_w) \hookrightarrow 
W(L^\infty,L^{p,q}(v)).
\]
Up to the authors knowledge this is a new convolution relation
on the $ax+b$-group even for $p,q\geq 1$.

\section*{Acknowledgement}

Major parts of this paper were developed during a stay 
at the Mathematical Institute of the University of Wroc{\l}aw, Poland. 
The author 
would like to thank its members for their warm hospitality. 
In particular, he would like to express his gratitude to M. Paluszy\'nski 
and R. Szwarc for interesting discussions. Many thanks go also to H.-Q. Bui for
hints on doubling weights and Muckenhoupt weights.
The author was supported 
by the European Union's Human Potential Programme under 
contracts HPRN-CT-2001-00273 (HARP) and  
HPRN-CT-2002-00285 (HASSIP).

%

\end{document}